\documentclass[12pt]{report}
\usepackage{amsmath, amssymb, amsthm, geometry, enumerate, graphicx, amsfonts, mathrsfs}
\theoremstyle{plain}

\theoremstyle{remark}

\numberwithin{equation}{section}

\begin{document}

\vspace{.2in}\parindent=0mm

\begin{flushleft}
 {\bf\Large {Construction of Nonuniform  Wavelet Frames
 \parindent=0mm \vspace{.1in}   on Non-Archimedean  Fields}}

\parindent=0mm \vspace{.3in}
  {\bf{Owais Ahmad$^{\star}$ and Neyaz Ahmad$^{\star\star}$} }
\end{flushleft}

\parindent=0mm \vspace{.1in}
{{\it\small$^{\star}${Department of Mathematics, National Institute of Technology, Srinagar-190006, Jammu and Kashmir, India. E-mail: $\text{siawoahmad@gmail.com}$}}

\parindent=0mm \vspace{.1in}
{{\it\small$^{\star\star}${Department of Mathematics, National Institute of Technology, Srinagar-190006, Jammu and Kashmir, India. E-mail: $\text{neyaznit@yahoo.co.in}$}}

\parindent=0mm \vspace{.2in}
{\bf{Abstract:}}
A constructive algorithm based on the theory of spectral pairs for constructing nonuniform wavelet basis in $L^2(\mathbb R)$ was considered by Gabardo and Nashed (J Funct. Anal. 158:209-241, 1998). In this setting, the associated translation set $\Lambda =\left\{ 0,r/N\right\}+2\,\mathbb Z$ is no longer a discrete subgroup of $\mathbb R$ but a spectrum associated with a certain one-dimensional spectral pair and the associated dilation is an even positive integer related to the given spectral pair. The main objective of this paper  is to develop  oblique and unitary extension principles  for the  construction  nonuniform wavelet frames over non-Archimedean Local fields of positive characteristic. An example  and some potential applications are also presented.

\parindent=0mm \vspace{.2in}

\parindent=0mm \vspace{.2in}

\parindent=0mm \vspace{.2in}
{\bf{1. Introduction}}

\parindent=0mm \vspace{.1in}
Duffin and Schaeffer \cite{duffCLAS} introduced the concept of frame in seperable Hilbert space while dealing with some deep problems in non-harmonic Fourier series. Frames are basis-like systems that span a vector space but allow for linear dependency, which can be used to reduce noise, find sparse representations, or obtain other desirable features unavailable with orthonormal bases. An important example about frame is wavelet frame, which is obtained by translating and dilating a finite family of functions. To mention only a few references on  wavelet frames, the reader is referred to \cite{{CK},{csLP},{chrisBOOK},{daubTEN},{DHRS},{dgmPNE}} and many references therein. Multiresolution analysis is an important mathematical tool since it provides a natural framework for understanding and constructing discrete wavelet systems.  The concept of MRA has been extended in various ways in recent years. These concepts are generalized to  $L^2\big(\mathbb R^d\big)$, to lattices different from  $\mathbb Z^d$, allowing the subspaces of MRA to be generated by Riesz basis instead of orthonormal basis, admitting a finite number of scaling functions, replacing the dilation factor 2 by an integer $M\geq 2$ or by an expansive matrix $A\in GL_{d}(\mathbb R)$ as long as $A\subset A\mathbb Z^d$. All these concepts are developed on regular lattices, that is the translation set is always a group. Recently, Gabardo and Nashed \cite{gabNON} considered a generalization of Mallat's \cite{malatMRA} celebrated theory of  MRA based on spectral pairs, in which the translation set acting on the scaling function associated with the  MRA to generate the subspace $V_0$ is no longer a group, but is the union of $\mathbb Z$ and a translate of $\mathbb Z$. Based on  one-dimensional spectral pairs, Gabardo and  Yu \cite{GN2} considered sets of   nonuniform wavelets  in $L^2(\mathbb{R})$.  

\parindent=8mm \vspace{.1in}
Ron and Shen \cite{Ron1} introduced the unitary extension principle which gives the construction of a multi-generated tight wavelet frame  for  $L^2(\mathbb{R}^d)$, based on a given refinable function. Tight wavelet frames gives more convenient way to represent a function in $L^2(\mathbb{R})$ in comparison of non-tight wavelet frames as in that case frame operator is constant multiple of identity operator in $L^2(\mathbb{R})$. Christensen and Goh   \cite{OCG1} generalized  the unitary extension principle to  locally compact abelian groups. They gave general constructions, based on B-splines on the group itself as well as on characteristic functions on the dual group. In real life application all signals are not obtained from  uniform shifts; so there is a natural question  regarding  analysis and decompositions of  this types of signals by a stable mathematical tool.   Gabardo and Nashed \cite{gabNON} and  Gabardo and  Yu \cite{GN2} filled this gap by  the concept of  nonuniform multiresolution analysis.  

\parindent=8mm \vspace{.1in}
A field  $\mathbb K$ equipped with a topology is called a local field if both the additive and multiplicative groups of $\mathbb K $  are locally compact abelian groups. The local fields are essentially classfied  into  two classes (excluding the connected local fields $\mathbb R $ and  $\mathbb C$). The local fields of characteristic zero include the $p$-adic field $\mathbb Q_p$. Examples of local fields of positive characteristic are the Cantor dyadic group and the Vilenkin $p$-groups. The Vilenkin $p$-groups are also called $p$-series fields.

\parindent=8mm \vspace{.1in}
During the last two decades, there is a substantial body of work that has been concerned with the construction of wavelets on  local fields.  Even though the structures and metrics of local fields of zero and positive characteristics are similar, their wavelet and MRA (multiresolution analysis) theory are quite different.  For example,  R. L. Benedetto and J. J. Benedetto \cite{benWAV} developed a wavelet theory for local fields and related groups. They did not develop the multiresolution analysis (MRA) approach, their method is based on the theory of wavelet sets and only allows the construction of wavelet functions whose Fourier transforms are characteristic functions of some sets. Khrennikov, Shelkovich and Skopina \cite{17} constructed a number of scaling functions generating an MRA of $L^2(\mathbb Q_p$). But later on in \cite{2}, Albeverio, Evdokimov and Skopina proved that all these scaling functions lead to the same Haar MRA and that there exist no other orthogonal test scaling functions generating an MRA except those described in \cite{17}. Some wavelet bases for $L^2(\mathbb Q_p$) different from the Haar system were constructed in \cite{{1},{9}} . These wavelet bases were obtained by relaxing the basis condition in the definition of an MRA and  form Riesz bases without any dual wavelet systems. For some related works on wavelets and frames on $\mathbb Q_p$, we refer to \cite{{3},{16},{18},{19}}. On the other hand, Lang \cite{{20},{21},{22}} constructed several examples of compactly supported wavelets for the Cantor dyadic group. Farkov  \cite{{10},{11}} has constructed many examples of wavelets for the Vilenkin $p$-groups. Jiang et al.\cite{jangMRA} pointed out a method for constructing orthogonal wavelets on  local field $\mathbb K$ with a constant generating sequence and derived necessary and sufficient conditions for a solution of the refinement equation to generate a multiresolution analysis of $L^2(\mathbb K)$. During the last two decades, p-adics has been extensively applied to a variety of problems in theoretical physics (string theory, cosmology, quantum theory, and disordered systems,) and biology (in modeling the thinking process and in genetics)\cite{{z8},{z9},{z7},{z1},{z2},{z3},{z4},{z5},{z6}}.

\parindent=8mm \vspace{.1in}
 Recently, Shah and Abdullah \cite{shahNUMRA} have generalized the concept of multiresolution analysis on Euclidean spaces $\mathbb R^n$ to  nonuniform multiresolution analysis on local fields of positive characteristic, in which the translation set acting on the scaling function associated with the multiresolution analysis to generate the subspace $V_{0}$ is no longer a group, but is the union of ${\mathcal Z}$ and a translate of ${\mathcal Z}$, where ${\mathcal Z}=\{u(n): n\in\mathbb N_{0}\}$  is a complete list of (distinct) coset representation of the unit disc $\mathfrak D$ in the locally compact Abelian group $\mathbb K^+.$  More precisely, this set is of the form $\Lambda=\left\{0,r/N \right\}+{\mathcal Z}$, where $N \ge 1$ is an integer and $r$ is an odd integer such that $r$ and $N$ are relatively prime. They call this a {\it  nonuniform multiresolution analysis} on local fields of positive characteristic.The notion of nonuniform wavelet frames on non-Archimedean local fields was introduced by Ahmad and Sheikh \cite{nuwf} and established a complete characterization of tight nonuniform wavelet frames on non-Archimedean local  fields. More results in this direction can also be found in \cite{{oAfr},{ogbr},{denNAS},{oJGP},{oSIS}} and the references therein. Drawing the inspiration from the above work, we develop the extension principles for the construction of nonuniform Parseval wavelet frames in $L^2(\mathbb K).$  

\parindent=8mm \vspace{.1in}
The remainder of the paper is as follows. In Section 2, we discuss preliminary results on non-Archimedean local fields. Section 3 is devoted to main results of this paper. Some potential applications are presented in Section 4.

\parindent=0mm \vspace{.2in}
{\bf{2. Preliminaries on  Non-Archimedean Local  Fields }}

\parindent=0mm \vspace{.1in}
{\it {2.1. Non-Archimedean Local Fields}}

\parindent=0mm \vspace{.1in}
A  non-Archimedean local field $\mathbb K$ is a  locally compact, non-discrete and totally disconnected field. If it is of characteristic zero, then  it is a field of $p$-adic numbers $\mathbb Q_p$ or its finite extension. If $\mathbb K$ is of positive characteristic, then $\mathbb K$ is a field of formal Laurent series over a finite field $GF(p^c)$. If $c =1$, it is a $p$-series field, while for $c\ne 1$, it is an algebraic extension of degree  $c$ of a $p$-series field. Let $\mathbb K$ be a fixed non-Archimedean local field  with the ring of integers ${\mathfrak D}= \left\{x \in K: |x| \le 1\right\}$. Since $K^{+}$ is a locally compact Abelian group, we choose a Haar measure $dx$ for $K^{+}$. The  field $K$ is locally compact, non-trivial, totally disconnected and complete topological field endowed with non--Archimedean norm  $|\cdot|:\mathbb K  \to \mathbb R^+$ satisfying

\parindent=0mm \vspace{.1in}
(a) $|x|=0$ if and only if $x = 0;$

\parindent=0mm \vspace{.1in}
(b) $|x\,y|=|x||y|$ for all $x, y\in \mathbb K$;

\parindent=0mm \vspace{.1in}
(c) $|x+y|\le \max \left\{ |x|, |y|\right\}$ for all $x, y\in \mathbb K$.

\parindent=0mm \vspace{.1in}
Property (c) is called the ultrametric inequality. Let ${\mathfrak B}= \left\{x \in \mathbb K: |x| < 1\right\}$ be the prime ideal of the ring of integers ${\mathfrak D}$ in $\mathbb K$. Then, the residue space ${\mathfrak D}/{\mathfrak B}$ is isomorphic to a finite field $GF(q)$, where $q = p^{c}$ for some prime $p$ and $c\in\mathbb N$. Since  $K$ is totally disconnected and $\mathfrak B$ is both prime and principal ideal, so there exist a prime element $\mathfrak p$ of $\mathbb K$ such that ${\mathfrak B}= \langle \mathfrak p \rangle=\mathfrak p {\mathfrak D}$. Let ${\mathfrak D}^*= {\mathfrak D}\setminus {\mathfrak B }=\left\{x\in \mathbb K: |x|=1   \right\}$. Clearly,  ${\mathfrak D}^*$ is a group of units in $\mathbb K^*$ and if $x\not=0$, then can write $x=\mathfrak p^n y, y\in {\mathfrak D}^*.$ Moreover, if ${\cal U}= \left\{a_m:m=0,1,\dots,q-1 \right\}$ denotes the fixed full set of coset representatives of ${\mathfrak B}$ in ${\mathfrak D}$, then every element $x\in K$ can be expressed uniquely  as $x=\sum_{\ell=k}^{\infty} c_\ell \,\mathfrak p^\ell $ with $c_\ell \in {\cal U}.$ Recall that ${\mathfrak B}$ is compact and open, so each  fractional ideal ${\mathfrak B}^k= \mathfrak p^k {\mathfrak D}=\left\{x \in K: |x| < q^{-k}\right\}$  is also compact and open and is a subgroup of $K^+$. We use the notation in Taibleson's book \cite{table}. In the rest of this paper, we use the symbols $\mathbb N, \mathbb N_0$ and $\mathbb Z$ to denote the sets of natural, non-negative integers and integers, respectively.

\parindent=8mm \vspace{.1in}
 Let $\chi$ be a fixed character on $K^+$ that is trivial on ${\mathfrak D}$ but  non-trivial on  ${\mathfrak B}^{-1}$. Therefore, $\chi$ is constant on cosets of ${\mathfrak D}$ so if $y \in {\mathfrak B}^k$, then $\chi_y(x)=\chi(y,x), x\in K.$ Suppose that $\chi_u$ is any character on $K^+$, then the restriction $\chi_u|{\mathfrak D}$ is a character on ${\mathfrak D}$. Moreover, as characters on ${\mathfrak D}, \chi_u=\chi_v$ if and only if $u-v\in {\mathfrak D}$. Hence, if  $\left\{u(n): n\in\mathbb N_0\right\}$ is a complete list of distinct coset representative of ${\mathfrak D}$ in $K^+$, then, as it was proved in \cite{table}, the set  $\left\{\chi_{u(n)}: n\in\mathbb N_0\right\}$   of distinct characters on ${\mathfrak D}$ is a complete orthonormal system on ${\mathfrak D}$.

\parindent=8mm \vspace{.1in}
We now impose a natural order on the sequence $\{u(n)\}_{n=0}^\infty$. We have ${\mathfrak D}/ \mathfrak B \cong GF(q) $ where $GF(q)$ is a $c$-dimensional vector space over the field $GF(p)$. We choose a set $\left\{1=\zeta_0,\zeta_1,\zeta_2,\dots,\zeta_{c-1}\right\}\subset {\mathfrak D^*}$ such that span$\left\{\zeta_j\right\}_{j=0}^{c-1}\cong GF(q)$. For $n \in \mathbb N_0$ satisfying
$$0\leq n<q,~~n=a_0+a_1p+\dots+a_{c-1}p^{c-1},~~0\leq a_k<p,~~\text{and}~k=0,1,\dots,c-1,$$

\parindent=0mm \vspace{.1in}
we define
$$u(n)=\left(a_0+a_1\zeta_1+\dots+a_{c-1}\zeta_{c-1}\right){\mathfrak p}^{-1}.\eqno(2.1)$$

\parindent=0mm \vspace{.1in}
Also, for $n=b_0+b_1q+b_2q^2+\dots+b_sq^s, ~n\in \mathbb N_{0},~0\leq b_k<q,k=0,1,2,\dots,s$, we set

$$u(n)=u(b_0)+u(b_1){\mathfrak p}^{-1}+\dots+u(b_s){\mathfrak p}^{-s}.\eqno(2.2)$$

\parindent=0mm \vspace{.1in}
This defines $u(n)$ for all $n\in \mathbb N_{0}$. In general, it is not true that $u(m + n)=u(m)+u(n)$. But, if $r,k\in\mathbb N_{0}\; \text{and}\;0\le s<q^k$, then $u(rq^k+s)=u(r){\mathfrak p}^{-k}+u(s).$ Further, it is also easy to verify that $u(n)=0$ if and only if $n=0$ and $\{u(\ell)+u(k):k \in \mathbb N_0\}=\{u(k):k \in \mathbb N_0\}$ for a fixed $\ell \in \mathbb N_0.$ Hereafter we use the notation $\chi_n=\chi_{u(n)}, \, n\ge 0$.

\parindent=8mm \vspace{.2in}
Let the local field $\mathbb K$ be of characteristic $p>0$ and $\zeta_0,\zeta_1,\zeta_2,\dots,\zeta_{c-1}$ be as above. We define a character $\chi$ on $K$ as follows:
$$\chi(\zeta_\mu {\mathfrak p}^{-j})= \left\{
\begin{array}{lcl}
\exp(2\pi i/p),&&\mu=0\;\text{and}\;j=1,\\
1,&&\mu=1,\dots,c-1\;\text{or}\;j \neq 1.
\end{array}
\right. \eqno(2.3)$$

\parindent=0mm \vspace{.1in}
{\it {2.2.  Fourier Transforms on Non-Archimedean Local Fields}}

\parindent=0mm \vspace{.2in}
 The Fourier transform of $f \in L^1(K)$ is denoted by $\hat f(\xi)$ and defined  by
\begin{align*}
{\mathcal F}\big\{f(x)\big\}=\hat f(\xi)=\int_K f(x)\overline{ \chi_\xi(x)}\,dx.\tag{2.4}
\end{align*}
It is noted that
$$\hat f(\xi)= \displaystyle \int_K f(x)\,\overline{ \chi_\xi(x)}dx= \displaystyle \int_K f(x)\chi(-\xi x)\,dx.$$

\parindent=0mm \vspace{.1in}
The properties of Fourier transforms on non-Archimedean local field $\mathbb K$ are much similar to those of on the classical field $\mathbb R$. In fact, the Fourier transform on non-Archimedean local  fields of positive characteristic have the following properties:
\begin{itemize}
  \item The map $f\to \hat f$ is a bounded linear transformation of $L^1(\mathbb K)$ into $L^\infty(\mathbb K)$, and $\big\|\hat f\big\|_{\infty}\le \big\|f\big\|_{1}$.
  \item If $f\in L^1(\mathbb K)$, then $\hat f$ is uniformly continuous.
  \item If $f\in L^1(\mathbb K)\cap L^2(\mathbb K)$, then $\big\|\hat f\big\|_{2}=\big\|f\big\|_{2}$.
\end{itemize}

\parindent=0mm \vspace{.1in}
The Fourier transform of a function $f\in L^2(\mathbb K)$ is defined by
\begin{align*}
\hat f(\xi)= \lim_{k\to \infty} \hat f_{k}(\xi)=\lim_{k\to \infty}\int_{|x|\le q^{k}} f(x)\overline{ \chi_\xi(x)}\,dx,\tag{2.5}
\end{align*}

\parindent=0mm \vspace{.0in}
where  $f_{k}=f\,\Phi_{-k}$ and $\Phi_{k}$ is  the characteristic function of ${\mathfrak B}^{k}$.  Furthermore, if $f\in L^2(\mathfrak D)$, then we define the Fourier coefficients of $f$ as
\begin{align*}
\hat f\big(u(n)\big)=\int_{\mathfrak D} f(x) \overline{ \chi_{u(n)}(x)}\,dx.\tag{2.6}
\end{align*}

\parindent=0mm \vspace{.0in}
The series $\sum_{n\in \mathbb N_{0}} \hat f\big( u(n)\big) \chi_{u(n)}(x)$ is called the Fourier series of $f$. From the standard $L^2$-theory for compact Abelian groups, we conclude that the Fourier series of $f$ converges to $f$ in $L^2(\mathfrak D)$ and Parseval's identity holds:
\begin{align*}
\big\|f\big\|^2_{2}=\int_{\mathfrak D}\big|f(x)\big|^2 dx= \sum_{n\in \mathbb N_{0}} \left| \hat f\big(u(n)\big)\right|^2.\tag{2.7}
\end{align*}

\parindent=0mm \vspace{.0in}
{\it {2.3. Uniform MRA on Non-Archimedean Local Fields}}

\parindent=8mm \vspace{.1in}
In order to able to define the concepts of uniform MRA and wavelets on non-Archimedean  local fields, we need analogous notions of translation and dilation. Since $\bigcup_{j\in\mathbb Z} \mathfrak p^{-j} {\mathfrak D}=\mathbb K$, we can regard $\mathfrak p^{-1}$ as the dilation and since $\left\{u(n):n\in\mathbb N_{0}\right\}$ is a complete list of distinct coset representatives of $\mathfrak D$ in $K$, the set ${\mathcal Z}=\left\{u(n):n\in\mathbb N_{0}\right\}$ can be treated as the translation set. Note that $\Lambda$ is a subgroup of $\mathbb K^{+}$ and unlike the standard wavelet theory on the real line, the translation set is not a group. Let us recall the definition of a uniform MRA on non-Archimedean local  fields of positive characteristic introduced by Jiang et al. in  \cite{jangMRA}.

\parindent=0mm \vspace{.1in}
{\bf{Definition 2.1.}} Let $\mathbb K$ be a non-Archimedean local field of positive characteristic $p>0$ and $ \mathfrak p$ be a prime element of $\mathbb K$. An MRA of $L^2(\mathbb K)$ is a  sequence of closed subspaces $\{V_j:j\in \mathbb Z\}$ of $L^2(\mathbb K)$ satisfying the following properties:

\parindent=0mm \vspace{.1in}
(a)\quad $V_j \subset V_{j+1}\; \text{for all}\; j \in \mathbb Z;$

\parindent=0mm \vspace{.1in}
(b)\quad $\bigcup_{j\in \mathbb Z}V_j\;\text{is dense in}\;L^2(\mathbb K);$

\parindent=0mm \vspace{.1in}
(c)\quad $\bigcap_{j\in \mathbb Z}V_j=\{0\};$

\parindent=0mm \vspace{.1in}
(d)\quad $f(x) \in V_j\; \text{if and only if}\;f({\mathfrak p}^{-1}x) \in V_{j+1}\; \text{for all}\; j \in \mathbb Z;$

\parindent=0mm \vspace{.1in}
(e) There exists a function $\phi \in V_0$, such that $\left\{\phi\big(x-u(k)\big): k\in \mathbb N_0\right\}$ forms an orthonormal basis for $V_0$.

\parindent=5mm \vspace{.1in}
According to the standard scheme for construction of MRA-based wavelets, for each $j$, we define a wavelet space $W_{j}$ as the orthogonal complement of $V_{j}$ in $V_{j+1}$, i.e., $V_{j+1}=V_{j}\oplus W_{j}, \, j\in\mathbb Z$, where $W_{j}\perp V_{j}, \, j\in\mathbb Z$. It is not difficult to see that
$$f(x)\in W_{j} \quad \text{if and only if}\quad f(\mathfrak p^{-1}x)\in W_{j+1},\quad j\in\mathbb Z.\eqno(2.7)$$

\parindent=0mm \vspace{.0in}
Moreover, they are mutually orthogonal, and we have the following orthogonal decompositions:
$$L^2(K)= \bigoplus_{j\in\mathbb Z} W_{j}=V_{0}\oplus \left(\bigoplus_{j\ge 0}W_{j}\right).\eqno(2.8)$$

\parindent=0mm \vspace{.0in}
As in the case of $\mathbb R^n$, we expect the existence of $q-1$ number of functions $\psi_{1}, \psi_{2},\dots, \psi_{q-1}$ to form a set of basic wavelets. In view of (2.7) and (2.8), it is clear that if $\left\{\psi_{1}, \psi_{2},\dots, \psi_{q-1}\right\}$ is a set of function such that the system $\left\{\psi_{\ell}\big(x-u(k)\big): 1\le \ell\le q-1, k\in\mathbb N_{0} \right\}$ forms an orthonormal basis for $W_{0}$, then  $\left\{q^{j/2}\psi_{\ell}(\mathfrak p ^{-j}x-u(k)\big): 1\le \ell\le q-1,j\in\mathbb Z, k\in\mathbb N_{0} \right\}$ forms an orthonormal basis for $L^2(K)$.

 \parindent=0mm \vspace{.1in}
{\it {2.4. Nonuniform MRA  on Non-Archimedean Local Fields}}

\parindent=8mm \vspace{.1in}
 For an integer $N \ge 1$ and an odd integer $r$ with $1\leq r \leq qN-1$ such that $r$ and $N$ are relatively prime, we define $$\Lambda = \left\{ 0, \dfrac{u(r)}{N}\right\}+{\mathcal Z}.$$

 \parindent=0mm \vspace{.1in}
where ${\mathcal Z}=\left\{ u(n): n\in \mathbb N_{0}\right\}$. It is easy to verify that $\Lambda$ is not a group on non-Archimedean local field $\mathbb K$, but is the union of ${\mathcal Z}$ and a translate of ${\mathcal Z}.$ Following is the definition of nonuniform multiresolution analysis (NUMRA) on non-Archimedean local fields of positive characteristic given by Shah and Abdullah \cite{shahNUMRA}.

\parindent=0mm \vspace{.1in}
{\bf{Definition 2.2.}} For an integer $N \ge 1$ and an odd integer $r$ with $1\leq r \leq qN-1$ such that $r$ and $N$ are relatively prime, an associated  NUMRA on non-Archimedean local field $\mathbb K$ of positive characteristic is a sequence of closed subspaces $\left\{V_j: j\in\mathbb Z\right\}$ of $L^2(\mathbb K)$ such that the following properties hold:

\parindent=0mm \vspace{.2in}
(a)\quad $V_j \subset V_{j+1}\; \text{for all}\; j \in \mathbb Z;$

\parindent=0mm \vspace{.1in}
(b)\quad $\bigcup_{j\in \mathbb Z}V_j\;\text{is dense in}\;L^2(\mathbb K);$

\parindent=0mm \vspace{.1in}
(c)\quad $\bigcap_{j\in \mathbb Z}V_j=\{0\};$

\parindent=0mm \vspace{.1in}
(d)\quad $f(\cdot) \in V_j\; \text{if and only if}\;f({\mathfrak p}^{-1}N\cdot) \in V_{j+1}\; \text{for all}\; j \in \mathbb Z;$

\parindent=0mm\vspace{.1in}

(e)~ There exists a function $\phi$ in $V_0$ such that $\left\{ \phi (\cdot- \lambda ): \lambda \in \Lambda\right\}$, is a complete orthonormal basis for $V_0$.

\parindent=8mm \vspace{.1in}
It is worth noticing that, when $N = 1$, one recovers from the definition above the  definition of an MRA on non-Archimedean local fields of positive characteristic $p>0$. When, $N > 1$, the dilation is induced by $\mathfrak p^{-1}N$ and $|\mathfrak p^{-1}|=q$ ensures that $qN\Lambda\subset {\mathcal Z} \subset \Lambda$.

\parindent=8mm \vspace{.2in}
For every $j\in\mathbb Z$, define $W_{j}$ to be the orthogonal complement of $V_{j}$ in $V_{j+1}$. Then we have
$$V_{j+1}=V_{j}\oplus W_{j}\quad \text{and}\quad W_{\ell}\perp W_{\ell^\prime}\quad \text{if}~\ell\ne \ell^\prime.\eqno(2.7) $$

\parindent=0mm \vspace{.1in}
It follows that for $j>J$,
$$V_{j}=V_{J}\oplus \bigoplus_{\ell=0}^{j-J-1}W_{j-\ell}\,,\eqno(2.8) $$

\parindent=0mm \vspace{.1in}
where all these subspaces are orthogonal. By virtue of condition (b) in the Definition 2.2, this implies
$$L^2(\mathbb K)=\bigoplus_{j\in\mathbb Z}W_{j},\eqno(2.9) $$

\parindent=0mm \vspace{.1in}
a decomposition of $L^2(\mathbb K)$ into mutually orthogonal subspaces.

\parindent=8mm \vspace{.2in}
As in the standard scheme, one expects the existence of $qN -1$ number of functions so that their translation by elements of $\Lambda$ and dilations by the integral powers of ${\mathfrak p^{-1}}N$ form an orthonormal basis for $L^2(\mathbb K)$.

\parindent=8mm \vspace{.1in}
Let $a$ and $b$ be any two fixed elements in $\mathbb K$. Then, for any  prime  $\mathfrak p$ and $m,n\in\mathbb N_{0}$,let $D_{\mathfrak p},T_{u(n) a}$ and $E_{u(m) b}$   be the unitary operators acting on $f\in L^2(\mathbb K)$ defined by :
\begin{align*}
& T_{u(n) a}f(x)=f\big(x-u(n) a\big), \qquad ~\text{(Translation)}\\
&E_{u(m) b}f(x)=\chi\big(u(m) b x\big) f(x),~\quad \text{(Modulation)}\\
&D_{\mathfrak p}f(x)=\sqrt {qN} f\left(\mathfrak p^{-1}N x\right), ~\quad\quad \text{(Dilation)}.
\end{align*}

\parindent=0mm \vspace{.1in}
Then for any $f\in L^2(K)$, the following results can easily be verified:
\begin{align*}
&{\mathscr F}\big\{ T_{u(n) a} f(x) \big\} = E_{-u(n) a}{\mathscr F}\big\{f(x)\big\},\\
&{\mathscr F}\big\{ E_{u(m) b} f(x) \big\} = T_{u(m) b}{\mathscr F}\big\{f(x)\big\},\\
&{\mathscr F}\big\{ D_{\mathfrak p^j} f(x) \big\}= D_{\mathfrak p^{-j}}{\mathscr F}\big\{f(x)\big\}\\
&D_{\mathfrak p^j}T_{u(n)a}=T_{(qN)^{-j}u(n)a}D_{\mathfrak p^j} .`
\end{align*}

{\bf{3. Main Results}}

\parindent=0mm \vspace{.1in}
We start this section with the following definition

\parindent=0mm \vspace{.1in}
{\bf{Definition 3.1.}} Let $\{\psi_1,\psi_2,\dots,\psi_{qN-1}\}$ be a family of non-zero functions  in  $L^2(\mathbb{K})$ .   The system
\begin{align*}
\left\{D_{\mathfrak p^j}T_{\lambda}\psi_{\ell} \right\}_{j\in \mathbb{Z},\lambda \in \Lambda \atop 1 \le \ell \le qN-1} = \bigcup_{\ell=1}^{qN-1}\Big\{(qN)^{\frac{j}{2}}\psi_{\ell}(\mathfrak p^{-1}N)^j x-\lambda\Big\}_{j\in \mathbb{Z},\lambda \in \Lambda}
\end{align*}
is called a  nonuniform wavelet frame  for $L^2(\mathbb{K})$,  if   there exist finite  positive   constants $A$ and $B$ such that
\begin{align*}
A\|f\|^2 \le \sum_{\ell=1}^{qN-1} \sum_{j\in \mathbb{Z}} \sum_{\lambda \in \Lambda}  |\langle f,{D_{\mathfrak p^j}}T_\lambda \psi_{\ell} \rangle|^2 \le B\|f\|^2 \ \text{for all} \ f \in L^2(\mathbb{K}).
\end{align*}
	
\parindent=8mm \vspace{.1in}
 In formulation of the unitary extension principle there is long list of assumption, instead of writing each assumption again and again, we state all assumptions at once and call it \emph{general setup}.

\parindent=0mm \vspace{.1in}
{\bf{General Setup}:} 

\parindent=0mm \vspace{.1in}
{\bf{Assumption 1}:} Let $\psi_0\in L^2(\mathbb{K})$  be such that
$$ \hat{\psi_0}(\mathfrak p^{-1}N\xi)=m_0(\xi)\hat{\psi_0}(\xi), \  m_0(\xi)\in L^{\infty}(\mathbb{K});~ {\textit Supp}~  {\hat{\psi_0}}(\xi)\subseteq q^2N \mathfrak{D};
	\lim \limits_{\xi \to 0^+}\hat{\psi_0}(\xi)=1.$$
{\bf{Assumption 2}:} Let $m_1,m_2,\dots,m_{qN-1} \in L^{\infty}(\mathbb{K})$,  and define $\psi_1,\psi_2,\dots,\psi_{qN-1} \in L^2(\mathbb{K})$ such that
\begin{align*}
\widehat{\psi}_{\ell}(\mathfrak p^{-1}N\xi)= m_{\ell} (\xi)\hat{\psi_0}(\xi), \ 1 \le  \ell \le  qN-1.
\end{align*}
Let  ${\cal M}(\xi)$ be  a $qN \times 1$ matrix given by
\begin{align*}
{\cal M}(\xi)=\begin{bmatrix}
m_0(\xi)\\
m_1(\xi)\\
\vdots\\
m_{qN-1}(\xi)
\end{bmatrix}_{qN \times 1}.
\end{align*}
Then, the  collection  $\{\psi_{\ell}, m_{\ell}\}_{\ell=0}^{qN-1}$  satistying the above assumptions is called  a {\it general setup}.

\parindent=8mm \vspace{.1in}
Our aim is to find conditions on ${\cal M}(\xi)$ such that the nonuniform system  $\left\{{D_{{\mathfrak p}^{j}}}T_{\lambda}\psi_{\ell}\right\}_{j\in \mathbb{Z},\lambda \in \Lambda \atop  1 \le \ell \le qN-1}$ constitutes a nonuniform Parseval frame for $L^2(\mathbb{K})$.

\parindent=0mm \vspace{.1in}
{\bf {Lemma 3.1. }} For any  $f \in L^1(\mathbb{K})$, the function $\mathcal{P}f (x)=\sum\limits_{k\in \mathbb{N}_0} f\left(x + Nu(k)\right)$ is well defined, $N$-periodic and belongs to $L^1(\mathfrak D)$.

\begin{proof}
It is clear that $\mathcal{P}f (x)=\sum\limits_{k\in \mathbb{N}_0} f(x + Nu(k))$ is $N$-periodic.
For any  $f \in L^1(\mathbb{R})$, we have
\begin{align*}
\int\limits_{\mathfrak{D}} \sum\limits_{k\in \mathbb{N}_0}|f(x +Nu(k))|d\,\gamma=\int\limits_{\mathbb{K}} |f(x)|\, dx <\infty.
\end{align*}
Thus, $\mathcal{P}f (x)$ is well defined a.e. on $\mathbb{K}$, and also  belongs to $L^1(\mathfrak{D})$.
\end{proof}

{\bf{Lemma 3.2.}} Assume that
\begin{enumerate}[$(i)$]
\item  $\psi_0 \in L^2\mathbb{(K)}$,  $\lim\limits_{\xi \to 0^+} \hat{\psi_0}(\xi) =1$  and  Supp $\hat{\psi_0} (\xi)$ $\subseteq \mathfrak{p}\mathfrak{D}$;
\item   $f \in L^2 \mathbb{(K)}$  such that $\hat{f} \in C_c{\mathbb{(K)}}$.
\end{enumerate}
  Then,  for any $\epsilon > 0 $ there exist $J \in \mathbb{Z}$ such that
\begin{align*}
(1-\epsilon)\|f\|^2 \le\sum\limits_{\lambda \in \Lambda} |\langle f,{D_{\mathfrak p^j}}T_\lambda \psi_0 \rangle|^2 \le(1+\epsilon) \|f\|^2 \ \text{for all} \  j \ge J.
\end{align*}

\begin{proof}
For any  $j \in \mathbb{Z}$, $({D_{\mathfrak p^j}}\hat{f})\bar{\hat{\psi_0}}\in L^1(\mathbb{K})$. By  invoking Lemma 3.1.,it is clear that the function  $ \mathcal{P}({D_{\mathfrak p^j}}\hat{f})\bar{\hat{\psi_0}}$ is well defined.
 Moreover, for $\xi \in \mathfrak{D}$, we have
\begin{align*}
\mathcal{P}({D_{\mathfrak p^j}}\hat{f})\bar{\hat{\psi_0}} &= \sum\limits_{k\in \mathbb{N}_0}(({D_{\mathfrak p^j}}\hat{f})\bar{\hat{\psi_0}})(\xi-Nu(k))\\\
&=\sum\limits_{k\in \mathbb{N}_0} ({D_{\mathfrak p^j}}\hat{f})(\xi-Nu(k)) \bar{\hat{\psi_0}}(\xi-Nu(k)).
\end{align*}
Thus,   $ \mathcal{P}(L^j\hat{f})\bar{\hat{\psi_0}}$ is bounded by finite linear combinations of translates of $\bar{\hat{\psi_0}}$	and  \break $ \mathcal{P}({D_{\mathfrak p^j}}\hat{f})\bar{\hat{\psi_0}} \in L^2(\mathfrak{D})$.

Since
\begin{align*}
\left\langle f,{D_{\mathfrak p^j}}T_\lambda \psi_0 \right\rangle& =\langle \widehat{f},\widehat{{D_{\mathfrak p^j}}T_\lambda \psi_0} \rangle\\\
& =\langle\widehat{f},{D_{\mathfrak p^{-j}}}E_{-\lambda}\hat{\psi_0}\rangle\\\
 &=\langle {D_{\mathfrak p^j}}\hat{f},E_{-\lambda}\hat{\psi_0} \rangle.
\end{align*}
Using  the fact  Supp $\hat{\psi}_0(\xi)\subseteq q^2N\mathfrak{D}$ , we obtain
$$\begin{array}{rcl}
\displaystyle\sum\limits_{\lambda \in \Lambda}	|\langle f,{D_{\mathfrak p^j}}T_\lambda \psi_0 \rangle|^2 &=&\displaystyle\sum\limits_{\lambda \in \Lambda}	|\langle {D_{\mathfrak p^j}}\hat{f},E_{-\lambda}\hat{\psi_0} \rangle  |^2 \notag\\\\
&=&\displaystyle\sum\limits_{\lambda\in \cal Z}	|\langle{D_{\mathfrak p^j}}\hat{f},E_{-\lambda}\hat{\psi_0} \rangle  |^2 +\sum\limits_{\lambda \in \left( \frac{r}{N}+\cal Z\right)}	|\langle {D_{\mathfrak p^j}}\hat{f},E_{-\lambda}\hat{\psi_0} \rangle  |^2 \notag\\\\
&=&\displaystyle\sum\limits_{m\in\mathbb{N}_0} \Big| \int\limits_{\mathfrak{D}}\mathcal{P}(({D_{\mathfrak p^j}}\hat{f})  \bar{\hat{\psi_0}} )(\xi) \chi(\mathfrak{p}{u(m) \xi}) \,d\xi  \Big|^2\notag\\\\
&&\qquad\qquad +\displaystyle\sum\limits_{m\in\mathbb{N}_0} \Big| \int\limits_{\mathfrak{D}}\mathcal{P}(({D_{\mathfrak p^j}}\hat{f})  \bar{\hat{\psi_0}} )(\xi) \chi\left(\left(\frac{r}{N} +\mathfrak{p} u(m)\right)\xi\right) \,d\xi  \Big|^2 \notag \\\\
&=&\displaystyle\sum\limits_{m\in\mathbb{N}_0} \Big| \int\limits_{qN\mathfrak{D}}(({D_{\mathfrak p^j}}\hat{f})  \bar{\hat{\psi_0}} ) \chi(\mathfrak{p}{u(m) \xi}) \,d\xi \Big|^2\notag\\\\
 &&\qquad\qquad +\displaystyle\sum\limits_{m\in\mathbb{N}_0} \Big| \int\limits_{qN\mathfrak{D}}(({D_{\mathfrak p^j}}\hat{f})  \bar{\hat{\psi_0}} )  \chi\left(\left(\frac{r}{N} +\mathfrak{p} u(m)\right)\xi\right) \,d\gamma  \Big|^2 \notag\\\\
&=&\displaystyle\sum\limits_{m\in\mathbb{N}_0} \Big| \int\limits_{qN\mathfrak{D}}(({D_{\mathfrak p^j}}\hat{f})  \bar{\hat{\psi_0}} ) \chi(\mathfrak{p}{u(m) \xi}) \,d\xi  \Big|^2\notag\\\\
  &&\qquad\qquad+\displaystyle\sum\limits_{m\in\mathbb{N}_0} \Big| \int\limits_{qN\mathfrak{D}}(({D_{\mathfrak p^j}}\hat{f})  \bar{\hat{\psi_0}} ) \chi\left(\left(\frac{r}{N} +\mathfrak{p} u(m)\right)\xi\right) \,d\xi  \Big|^2.
\end{array}$$ 
By using the Parseval identity  on $L^2(\mathfrak{pD})$ with respect to an orthonormal bases
$\left\{\sqrt{q}\;\chi\left(\mathfrak{p}{u(m) \xi}\right)\right\}$  in (3.1), we obtain
\begin{align}
\sum\limits_{\lambda \in \Lambda}	|\langle f,{D_{\mathfrak p^j}}T_\lambda \psi_0 \rangle|^2
&=\frac{1}{2}\int\limits_{\mathfrak{pD}}|({D_{\mathfrak p^j}}\hat{f})  \bar{\hat{\psi_0}}|^2\,d\xi +\frac{1}{2}\int\limits_{\mathfrak{pD}}|({D_{\mathfrak p^j}}\hat{f})  \bar{\hat{\psi_0}}|^2\,d\xi	\notag\\
&=\int\limits_{\mathfrak{pD}}|({D_{\mathfrak p^j}}\hat{f})  \bar{\hat{\psi_0}}|^2\,d\xi.\tag{3.2}
\end{align}

For a given  $\epsilon >0$ , we can choose  $\gamma \in \mathfrak{pB} $ in such a manner that
$$(1-\epsilon)\le |\hat{\psi_0}(\xi)|^2 \le (1+\epsilon),  \  \text{where} \  0<\xi<\gamma.\eqno(3.3)$$

Choose $J \in \mathbb{Z}$ large enough,  so that Supp $({D_{\mathfrak p^j}}\hat{f}) \subset \gamma\mathfrak{D}$ for all $j \geq J$. Then, by (3.2), we have
$$\sum\limits_{\lambda \in \Lambda}	|\langle f,{D_{\mathfrak p^j}}T_\lambda \psi_0 \rangle|^2=
\int\limits_{\gamma\mathfrak{D}}|({D_{\mathfrak p^j}}\hat{f})  \bar{\hat{\psi_0}}|^2\,d\xi \  \text{for all} \ j \geq J.\eqno(3.4)$$

 Using the fact that ${D_{\mathfrak p^j}}$ is unitary map and the equations (3.3), (3.4), we have
\begin{align*}
(1-\epsilon)\|\hat{f}\|^2 \le \sum\limits_{\lambda \in \Lambda}	|\langle f,{D_{\mathfrak p^j}}T_\lambda \psi_0 \rangle|^2 \le (1+\epsilon)\|\hat{f}\|^2 \ \text{ for all} \  j \ge J.
\end{align*}
 Since the  Fourier transform is  unitary map, we get
\begin{align*}
(1-\epsilon)\|{f}\|^2 \le \sum\limits_{\lambda \in \Lambda}	|\langle f,{D_{\mathfrak p^j}}T_\lambda \psi_0 \rangle|^2 \le (1+\epsilon)\|{f}\|^2 \ \text{ for all} \  j \ge J.
\end{align*}
This completes  the proof of the lemma 3.2.
\end{proof}

{\bf{Lemma 3.3.}} Suppose that
\begin{enumerate}[$(i)$]
\item  $\psi_0 \in L^2(\mathbb{K})$ satisfies  Supp $\hat{\psi_0} \subseteq q^2N\mathfrak{D}$ and
$\hat{\psi_0}(\mathfrak{p}^{-1}N\xi)=m_0(\xi)\hat{\psi_0}(\xi)$,where $m_0(\xi) \in L^\infty (\mathbb{K})$;
\item  $f \in L^2(\mathbb{K})$ with $\hat{f} \in C_c(\mathbb{K})$,  and $ m_1,m_2,\dots,m_{qN-1} \in  L^\infty (\mathbb{K})$  such that  the $qN\times 1$ matrix
$${\cal M}(\xi)=
\begin{bmatrix}
m_0(\xi)\\
m_1(\xi)\\
\vdots\\
m_{qN-1}(\xi)
\end{bmatrix}_{qN\times 1}$$
satisfies ${{\cal M}(\xi)}^\ast {\cal M}(\xi)=1$ $a.e.$;
\item $\psi_1,\psi_2, \dots \psi_{qN-1} \in L^2(\mathbb{K})$ such that $\hat{\psi}_{\ell}(\mathfrak{p}^{-1} N\xi)=m_{\ell}(\xi)\hat{\psi_0}(\xi),  \  1 \le \ell \le qN-1$.
\end{enumerate}
Then
\begin{align*}
\sum \limits_{\ell=0}^{qN-1}\sum\limits_{\lambda \in \Lambda} |\langle f,{D_{\mathfrak p^{j-1}}}T_\lambda \psi_{\ell} \rangle|^2 =\sum\limits_{\lambda \in \Lambda} |\langle f,{D_{\mathfrak p^j}}T_\lambda \psi_0 \rangle|^2.
\end{align*}

\begin{proof}
For any  $ j\in \mathbb{Z}$ and for any  $ 0 \le \ell \le qN-1$, we  have
\begin{align}
 \langle f,{D_{\mathfrak p^{j-1}}}T_\lambda \psi_{\ell} \rangle &={D_{\mathfrak p^{-j}}}f,{D_{\mathfrak p^{-1}}}T_\lambda \psi_{\ell} \rangle \notag\\\
&=\langle {D_{\mathfrak p^{-j}}}f,T_{(qN)\lambda}{D_{\mathfrak p^{-1}}}\psi_{\ell} \rangle \notag\\
&=\langle {D_{\mathfrak p^{-j}}}\hat{f},E_{-(qN)\lambda}D\hat\psi_\ell\rangle \notag\\
&=\int \limits_{\mathbb{K}}({D_{\mathfrak p^j}}\hat{f})(\gamma)\sqrt{qN}\overline{\hat{\psi_{\ell}}(\mathfrak{p}^{-1} N\xi)}\;\chi\left(\mathfrak{p}^{-1}N \lambda\xi\right) d\xi \notag\\
&=\sqrt{qN}\int \limits_{\mathbb{K}}({D_{\mathfrak p^j}}\hat{f})(\xi)\overline{m_{\ell}(\xi)}\overline{\hat{\psi_0}(\xi)}\;\chi\left(\mathfrak{p}^{-1}N \lambda\xi\right) d\xi. \tag{3.5}
\end{align}
\newpage
Apply  Parseval identity  on $ L^2(q^2N\mathfrak{D} )$ with respect to orthonormal basis\break $\left\{q\sqrt{N} \chi\left( \mathfrak{p}^{-2} Nu(m) \xi\right); m \in \mathbb N_0\right\}$, and the fact Supp $\hat{\psi_0}\subseteq q^2 N\mathfrak{D}$, we have

\begin{align}
\sum\limits_{\lambda \in \Lambda} |\langle f,{D_{\mathfrak p^{j-1}}}T_\lambda \psi_{\ell} \rangle|^2 &=\sum\limits_{\lambda \in \mathcal Z} |\langle f,{D_{\mathfrak p^{j-1}}}T_\lambda \psi_{\ell} \rangle|^2+\sum\limits_{\lambda \in (\frac{r}{N}+\mathcal{Z})} |\langle f,{D_{\mathfrak p^{j-1}}}T_\lambda \psi_{\ell} \rangle|^2 \notag\\
&=\sum\limits_{m\in\mathbb{Z}} \Big|\sqrt{qN} \int\limits_{\mathfrak{D}}\mathcal{P}(({D_{\mathfrak p^j}}\hat{f})\overline{m_{\ell}(\xi)}\overline{(\hat{\psi_0}})(\xi)\; \chi\left( \mathfrak{p}^{-2} N u(m) \xi\right)\, d\xi  \Big|^2 \notag \\
& \quad +\sum\limits_{m\in\mathbb{N}_0} \Big|\sqrt{qN} \int\limits_{\mathfrak{D}}\mathcal{P}(({D_{\mathfrak p^j}}\hat{f})\overline{m_{\ell}(\xi)}\overline{(\hat{\psi_0}})(\xi)\; \chi\left( \mathfrak{p}^{-2} N \left(\frac{r}{N}+u(m)\right) \xi\right)\, d\xi \Big|^2  \notag\\
&=\frac{1}{2}\sum\limits_{m\in\mathbb{N}_0} \Big| \int\limits_{q^2N\mathfrak{D}}({D_{\mathfrak p^j}}\hat{f})(\xi)\overline{m_{\ell}(\xi)}\overline{\hat{\psi_0}(\xi)}\;\chi\left(  \mathfrak{p}^{-2}Nu(m) \xi\right) q \sqrt{N}\, d\xi \Big|^2 \notag\\
& \quad +\frac{1}{2}\sum\limits_{m\in\mathbb{N}_0} \Big| \int\limits_{q^2N\mathfrak{D}}({D_{\mathfrak p^j}}\hat{f})(\xi)\overline{m_{\ell}(\xi)}\overline{\hat{\psi_0}(\xi)}\chi\left( \left(\mathfrak{p^{-1}}u(r)+\mathfrak{p}^{-2}Nu(m)\right)\xi \right) q \sqrt{N}\, d\xi  \Big|^2 \notag\\
&=\frac{1}{2}\int\limits_{q^2N\mathfrak{D}}|({D_{\mathfrak p^j}}\hat f)(\xi)\overline{m_\ell(\xi)}\overline{\hat{\psi_0}(\xi)}|^2 d\xi+\frac{1}{2}\int\limits_{q^2N\mathfrak{D}}|({D_{\mathfrak p^j}}\hat f)(\xi)\overline{m_\ell(\xi)}\overline{\hat{\psi_0}(\xi)}|^2 d\xi \notag\\
&=\int\limits_{q^2N\mathfrak{D}}|({D_{\mathfrak p^j}}\hat f)(\xi)\overline{m_\ell(\xi)}\overline{\hat{\psi_0}(\xi)}|^2 d\xi.\tag{3.6}
\end{align}

Since ${\cal M}(\xi)^\ast {\cal M}(\xi)=1$ a.e., so ${\cal M}(\xi)$ could be consider as an isometry from $\mathbb{C}^{1}$ into $\mathbb{C}^{qN}$. Using (3.6), we have
\begin{align}
\sum \limits_{\ell=0}^{qN-1}\sum\limits_{\lambda \in \Lambda} |\langle f,{D_{\mathfrak p^{j-1}}}T_\lambda \psi_{\ell} \rangle|^2 & =\sum \limits_{\ell=0}^{qN-1}   \int\limits_{q^2N\mathfrak{D}}|({D_{\mathfrak p^j}}\hat f)(\xi)\overline{m_\ell(\xi)}\overline{\hat{\psi_0}(\xi)}|^2 d\xi\notag\\
&=\displaystyle{\int \limits_{q^2N\mathfrak{D}}\left\|\begin{bmatrix}
\overline{m_0(\xi)}\notag\\
\vdots\\
\overline{m_{qN-1}(\xi)}
\end{bmatrix}_{qN \times 1} [({D_{\mathfrak p^j}}\hat{f})\overline{\hat{\psi_0}}]_{1\times 1}\right\|^2_{\mathbb{C}^{qN}}}d\xi \notag\\
&=\int\limits_{q^2N\mathfrak{D}}\left\|\overline{{\cal M}(\xi)}_{qN \times1}[({D_{\mathfrak p^j}}\hat{f})\overline{\hat{\psi_0}}]_{1\times 1}\right\|^2_{\mathbb{C}^{qN}}d\xi \notag\\
&=\int \limits_{q^2N\mathfrak{D}}|({D_{\mathfrak p^j}}\hat{f})(\xi)\overline{\hat{\psi_0}}(\xi)|^2\, d\xi.\tag{3.7}
\end{align}
Also
\begin{align}\label{eq3.99xv}
&\sum\limits_{\lambda \in \Lambda} |\langle f,{D_{\mathfrak p^j}}T_\lambda \psi_0 \rangle|^2 =\sum\limits_{\lambda \in 2\mathbb{Z}} |\langle f,{D_{\mathfrak p^j}}T_\lambda \psi_0 \rangle|^2+\sum\limits_{\lambda \in (\frac{r}{N}+2\mathbb{Z})} |\langle f,{D_{\mathfrak p^j}}T_\lambda \psi_0 \rangle|^2 \notag\\
&\qquad\qquad\qquad\qquad=\sum\limits_{m\in\mathbb{N}_0} \Big| \int\limits_{\mathbb{K}}({D_{\mathfrak p^j}}\hat{f})(\xi)\overline{\hat{\psi_0}(\xi)}\; \chi\left(\mathfrak{p}^{-1}u(m)\xi\right)\, d\xi  \Big|^2\notag\\
 &\qquad\qquad\qquad\qquad\qquad\qquad+\sum\limits_{m\in\mathbb{N}_0} \Big| \int\limits_{\mathbb{K}}({D_{\mathfrak p^j}}\hat{f})(\xi)\overline{\hat{\psi_0}(\xi)}\; \chi\left(\left(\frac{r}{N}+\mathfrak{p}^{-1} u(m)\right)\right) \, d\xi  \Big|^2.\tag{3.8}
\end{align}
 Applying the  Parseval formula on $L^2(\mathfrak{pD})$ with respect to orthonormal basis $\{\sqrt{q}\;\chi\left(\mathfrak{p}^{-1}u(m)\xi\right)\}_{m\in \mathbb{Z}}$  and the fact  Supp $\hat{\psi_0} \subseteq q^2N\mathfrak{D}$,  we obtain
$$\begin{array}{rcl}
\displaystyle\sum\limits_{\lambda \in \Lambda} |\langle f,{D_{\mathfrak p^j}}T_\lambda \psi_0 \rangle|^2 &=&\displaystyle\sum\limits_{m\in\mathbb{N}_0} \Big| \int\limits_{\mathfrak{D}}\mathcal{P}(({D_{\mathfrak p^j}}\hat{f})\overline{\hat{\psi_0})}(\xi)\chi\left(\mathfrak{p}^{-1}u(m)\xi\right)\, d\xi  \Big|^2\notag\\\
&&\qquad\qquad +\displaystyle\sum\limits_{m\in\mathbb{N}_0} \Big| \int\limits_{\mathfrak{D}}\mathcal{P}(({D_{\mathfrak p^j}}\hat{f})\overline{\hat{\psi_0})}(\xi)\chi\left(\left(\frac{r}{N}+\mathfrak{p}^{-1} u(m)\right)\right) \, d\xi  \Big|^2 \quad \Big(\text{by} ~(3.8)\Big) \notag\\\
&=&\dfrac{1}{2}\displaystyle\sum\limits_{m\in\mathbb{Z}} \Big| \int\limits_{0}^{\frac{1}{2}}({D_{\mathfrak p^j}}\hat{f})(\xi)\overline{\hat{\psi_0}}(\xi)\sqrt{q}\chi\left(\mathfrak{p}^{-1}u(m)\xi\right)\, d\xi  \Big|^2\notag\\\
&&\qquad\qquad +\dfrac{1}{2} \displaystyle\sum\limits_{m\in\mathbb{N}_0} \Big| \int\limits_{\mathfrak{pD}}({D_{\mathfrak p^j}}\hat{f})(\xi)\overline{\hat{\psi_0}}(\xi)\sqrt{q}\chi\left(\left(\frac{r}{N}+\mathfrak{p}^{-1} u(m)\right)\right) \, d\xi  \Big|^2 \notag\\\\
&=&\dfrac{1}{2}\displaystyle\int\limits_{\mathfrak{pD}}|({D_{\mathfrak p^j}}\hat f)(\xi)\overline{\hat{\psi_0}(\xi)} |^2d\xi+\frac{1}{2}\displaystyle\int\limits_{\mathfrak{pD}}|({D_{\mathfrak p^j}}\hat f)(\xi)\overline{\hat{\psi_0}(\xi)} |^2d\xi \notag\\\\
&=&\displaystyle\int\limits_{\mathfrak{pD}}|({D_{\mathfrak p^j}}\hat f)(\xi)\overline{\hat{\psi_0}(\xi)} |^2d\xi\notag\\\\
&=&\displaystyle\int\limits_{q^2N\mathfrak{D}}|({D_{\mathfrak p^j}}\hat f)(\xi)\overline{\hat{\psi_0}(\xi)} |^2d\xi. \qquad\qquad\qquad\qquad\quad\qquad\qquad\qquad{(3.9)}
\end{array}$$
The proof clearly follows from   (3.7) and (3.9).
\end{proof}

{\bf{Lemma 3.4.}} Let $\{\psi_{\ell}, m_{\ell}\}_{\ell=0}^{qN-1}$ be a general setup, and let  ${\cal M}(\xi)^{\ast}{\cal M}(\xi)$=1. Then,  the following holds.
	\begin{enumerate}[$(i)$]
\item $\{T_\lambda \psi_0\}_{\lambda \in \Lambda} $ is Bessel sequence with Bessel bound $1$. \label{2.3i}
\item For any $f \in L^2(\mathbb{K})$,\label{2.3ii}
\begin{align*}
\lim \limits_{j \to -\infty}\sum\limits_{\lambda \in \Lambda} |\langle f,{D_{\mathfrak p^j}}T_\lambda \psi_0 \rangle|^2=0.
\end{align*}
	\end{enumerate}

\begin{proof} $(i):$  Let  $f \in L^2(\mathbb{K}) $ be  such that $\hat{f} \in C_c(\mathbb{K}) $. For a given  $\epsilon >0$ ,  by Lemma 3.2, we can find an integer $j>0$ such that

$$\sum\limits_{\lambda \in \Lambda} |\langle f,{D_{\mathfrak p^j}}T_\lambda \psi_0 \rangle|^2 \le (1+\epsilon)\|f\|^2.\eqno(3.10)$$

Also, by Lemma 3.2, we have
$$\sum\limits_{\lambda \in \Lambda} |\langle f,{D_{\mathfrak p^{j-1}}}T_\lambda \psi_{\ell} \rangle|^2 \le \sum\limits_{\lambda \in \Lambda} |\langle f,{D_{\mathfrak p^j}}T_\lambda \psi_0 \rangle|^2.\eqno(3.11)$$
 On applying (3.11) $j$ times and  using (3.10), we get
\begin{align*}
\sum\limits_{\lambda \in \Lambda} |\langle f,T_\lambda \psi_0 \rangle|^2 \le
		\sum\limits_{\lambda \in \Lambda} |\langle f,{D_{\mathfrak p^j}}T_\lambda \psi_0 \rangle|^2 \le (1+\epsilon)\|f\|^2 .
\end{align*}
  Since $\epsilon >0$ was arbitrary, we have
\begin{align*}
\sum\limits_{\lambda \in \Lambda} |\langle f,T_\lambda \psi_0 \rangle|^2 \le\|f\|^2.
	\end{align*}
	
Since this inequality holds on a dense subset of $L^2(\mathbb{K})$,therefore  it holds on $L^2(\mathbb{K})$. This proves (i) of the lemma.

\parindent=0mm \vspace{.1in}
$(ii)$:  Let $f\in L^2(\mathbb{K})$. Since ${D_{\mathfrak p^j}}$ is unitary map  for all $j \in \mathbb{Z}$, by using (i), the family  $\{{D_{\mathfrak p^j}}T_\lambda \psi_0\}_{\lambda \in \Lambda} $ is Bessel sequence with Bessel bound $1$. For any  $j \in \mathbb{Z}$ and  for any bounded set  $\Delta \subset \mathbb{K}$, we have
\begin{align*}
\sum\limits_{\lambda \in \Lambda} |\langle f,{D_{\mathfrak p^j}}T_\lambda \psi_0 \rangle|^2 &\le 
\sum\limits_{\lambda \in \Lambda} |\langle f{\bf 1}_{\Delta},{D_{\mathfrak p^j}}T_\lambda \psi_0 \rangle|^2 +  \sum\limits_{\lambda \in \Lambda} |\langle f(1-{\bf 1}_{\Delta}),{D_{\mathfrak p^j}}T_\lambda \psi_0 \rangle|^2\\
&\le 
\sum\limits_{\lambda \in \Lambda} |\langle {\bf 1}_{\Delta},{D_{\mathfrak p^j}}T_\lambda \psi_0 \rangle|^2+\|f(1-{\bf 1}_{\Delta})\|^2.
\end{align*}
Now, $\|f(1-{\bf 1}_{\Delta})\|^2 \to 0$,  if we choose $\Delta$ to be sufficiently large. Therefore, we only need to show
\begin{align*}
\sum\limits_{\lambda \in \Lambda} |\langle f{\bf 1}_{\Delta},{D_{\mathfrak p^j}}T_\lambda \psi_0 \rangle|^2 \rightarrow 0 \ \text{as} \  j  \rightarrow  -\infty.
\end{align*}
Apply Cauchy-Schwarz's inequality for integrals, we get
\begin{align}
\sum\limits_{\lambda \in \Lambda} |\langle f{\bf 1}_{\Delta},{D_{\mathfrak p^j}}T_\lambda \psi_0 \rangle|^2 &=(qN)^j \sum\limits_{\lambda \in \Lambda}\Big|\int \limits_{\Delta} f(x)\overline{\psi_0((\mathfrak{p}^{-1}N)^j x-\lambda)} \,dx \Big|^2 \notag\\
&\le (qN)^j\|f\|^2 \sum \limits_{\lambda \in \Lambda} \int \limits_{\Delta} |\psi_0((\mathfrak{p}^{-1}N)^j x-\lambda)|^2 \,dx \notag\\
&=\|f\|^2 \sum \limits_{\lambda \in \Lambda} \displaystyle{\int \limits_{(qN)^j \Delta-\lambda} |\psi_0(x)|^2 \,dx}.\tag{3.12}
\end{align}
On applying the Lebesgue  dominated convergence theorem in (3.12), we have
\begin{align*}
\sum\limits_{\lambda \in \Lambda} |\langle f{\bf 1}_{\Delta},{D_{\mathfrak p^j}}T_\lambda \psi_0 \rangle|^2 \rightarrow 0\  \text{as}\   j \to -\infty.
\end{align*}
Hence (ii) is proved.
\end{proof}

\parindent=0mm \vspace{.1in}
{\bf{Theorem 3.1.}} Let $\{\psi_{\ell}, m_{\ell}\}_{\ell=0}^{qN-1}$ be  a general setup and    ${\cal M}(\xi)^{\ast}{\cal M}(\xi)=1$. Then the   system $\{{D_{\mathfrak p^j}}T_{\lambda}\psi_{\ell}\}_{j\in \mathbb{Z},\lambda \in \Lambda \atop  1 \le \ell \le qN-1}$ constitutes a Parseval frame for $L^2(\mathbb{K})$.

\begin{proof}
Let $\epsilon >0$ be given. Consider a function $f \in L^2(\mathbb{K})$ such that $\hat{f} \in C_c{(\mathbb{K})}$. By Lemma 3.2., we can choose $J>0$ such that for all $j\geq J$,
$$(1-\epsilon)\|f\|^2 \le\sum\limits_{\lambda \in \Lambda} |\langle f,{D_{\mathfrak p^j}}T_\lambda \psi_0 \rangle|^2 \le(1+\epsilon) \|f\|^2.\eqno(3.13)$$
Using  Lemma 3.3., we have
$$\begin{array}{rcl}
\displaystyle\sum\limits_{\lambda \in \Lambda} |\langle f,{D_{\mathfrak p^j}}T_\lambda \psi_0 \rangle|^2 &=&\displaystyle\sum \limits_{\ell=0}^{qN-1}\sum\limits_{\lambda \in \Lambda} |\langle f,{D_{\mathfrak p^{j-1}}}T_\lambda \psi_{\ell} \rangle|^2 \notag\\
&=&\displaystyle\sum\limits_{\lambda \in \Lambda} |\langle f,{D_{\mathfrak p^{j-1}}}T_\lambda \psi_0 \rangle|^2+\sum \limits_{\ell=1}^{qN-1}\sum\limits_{\lambda \in \Lambda} |\langle f,{D_{\mathfrak p^{j-1}}}T_\lambda \psi_{\ell} \rangle|^2.\qquad\qquad{(3.14)}
\end{array}$$
Applying  Lemma 3.3., on $\sum\limits_{\lambda \in \Lambda} |\langle f,{D_{\mathfrak p^{j-1}}}T_\lambda \psi_0 \rangle|^2$, we get
$$\sum\limits_{\lambda \in \Lambda} |\langle f,{D_{\mathfrak p^{j-1}}}T_\lambda \psi_0 \rangle|^2=\sum\limits_{\lambda \in \Lambda} |\langle f,{D_{\mathfrak p^{j-2}}}T_\lambda \psi_0 \rangle|^2+\sum \limits_{\ell=1}^{qN-1}\sum\limits_{\lambda \in \Lambda} |\langle f,{D_{\mathfrak p^{j-2}}}T_\lambda \psi_{\ell} \rangle|^2.\eqno(3.15)$$
By using (3.14) and (3.15), it follows that
\begin{align*}
\sum\limits_{\lambda \in \Lambda} |\langle f,{D_{\mathfrak p^j}}T_\lambda \psi_0 \rangle|^2=\sum\limits_{\lambda \in \Lambda} |\langle f,{D_{\mathfrak p^{j-2}}}T_\lambda \psi_0 \rangle|^2+\sum\limits_{\ell=1}^{qN-1}\sum\limits_{\lambda \in \Lambda}\sum\limits_{\gamma=j-2}^{j-1}| \langle f,{D_{\mathfrak p^\gamma}}T_\lambda \psi_{\ell} \rangle|^2.
\end{align*}
On repeating the above arguments, for any $m <j$, we have
$$\sum\limits_{\lambda \in \Lambda} |\langle f,{D_{\mathfrak p^j}}T_\lambda \psi_0 \rangle|^2=\sum\limits_{\lambda \in \Lambda} |\langle f,{D_{\mathfrak p^m}}T_\lambda \psi_0 \rangle|^2+\sum\limits_{\ell=1}^{qN-1}\sum\limits_{\lambda \in \Lambda}\sum\limits_{\nu=m}^{j-1}| \langle f,{D_{\mathfrak p^\nu}}T_\lambda \psi_{\ell} \rangle|^2.\eqno(3.16)$$
It follows from (3.13) and (3.16) that  for all $j\ge J,m<j$,
\begin{align*}
(1-\epsilon)\|f\|^2 \le \sum\limits_{\lambda \in \Lambda} |\langle f,{D_{\mathfrak p^m}}T_\lambda \psi_0 \rangle|^2+\sum\limits_{\ell=1}^{qN-1}\sum\limits_{\lambda \in \Lambda}\sum\limits_{\nu=m}^{j-1}| \langle f,{D_{\mathfrak p^\nu}}T_\lambda \psi_{\ell} \rangle|^2 \le (1+\epsilon)\|f\|^2.
\end{align*}
Letting $m \to -\infty$ in above and using Lemma 3.4~(ii), we have
\begin{align}
(1-\epsilon)\|f\|^2 \le \sum\limits_{\ell=1}^{qN-1}\sum\limits_{\lambda \in \Lambda}\sum\limits_{\nu=-\infty}^{j-1}| \langle f,{D_{\mathfrak p^\nu}}T_\lambda \psi_{\ell} \rangle|^2 \le (1+\epsilon)\|f\|^2.\tag{3.17}
\end{align}
Letting $j \to \infty$  in (3.17), we have
\begin{align*}
(1-\epsilon)\|f\|^2 \le \sum\limits_{\ell=1}^{qN-1}\sum\limits_{\lambda \in \Lambda}\sum\limits_{\nu=-\infty}^{\infty}| \langle f,{D_{\mathfrak p^\nu}}T_\lambda \psi_{\ell} \rangle|^2 \le (1+\epsilon)\|f\|^2.
\end{align*}
Since $\epsilon >0$ was arbitrary, we obtain
\begin{align*}
\sum\limits_{\ell=1}^{qN-1}\sum\limits_{\lambda \in \Lambda}\sum\limits_{\nu \in \mathbb{Z}}| \langle f,{D_{\mathfrak p^\nu}} T_\lambda \psi_{\ell} \rangle|^2=\|f\|^2  \text{for all} \ f \in L^{2}(\mathbb{K}),
\end{align*}
as desired.
\end{proof}
\parindent=8mm \vspace{.1in}
Now we proceed to state and prove  oblique extension principle  for  the construction of nonuniform wavelet frames in $L^{2}(\mathbb{K})$.

\parindent=0mm \vspace{.1in}
{\bf{Theorem 3.2.}}   Let $\{\psi_{\ell}, m_{\ell}\}_{\ell=0}^{qN-1}$ be a general setup. Assume that  there exist strictly positive function $\varphi \in L^{\infty}(\mathbb{K})$ for which $ \lim\limits_{\xi \to 0^{+}} \varphi(\xi)=1,$ and 
$$\varphi(\mathfrak{p}^{-1} N \xi)|m_0(\xi)|^2+\sum\limits_{\ell=1}^{qN-1}|m_{\ell}(\xi)|^2=\varphi(\xi).$$
 Then the system $\{{D_{\mathfrak p^j}}T_{\lambda}\psi_{\ell}\}_{j\in \mathbb{Z},\lambda \in \Lambda \atop  1 \le \ell \le qN-1}$ is  a Parseval nonuniform wavelet frame for $L^2(\mathbb{K})$.

\begin{proof}
Define $\widetilde{\psi}_0\in L^2(\mathbb{K})$ such that
\begin{align}
\widehat{\widetilde{\psi}}_0(\xi)=\sqrt{\varphi (\xi)}\hat{\psi_0} (\xi).\tag{3.18}
\end{align}
Define functions $\widetilde{m}_0,\widetilde{m}_1,\dots,\widetilde{m}_{qN-1}$ as follows
\begin{align*}
\widetilde{m_0}(\xi)&=\sqrt{\frac{\varphi(\mathfrak{p}^{-1}N\xi)}{\varphi(\xi)}}m_0(\xi),\\
\widetilde{m_{\ell}}(\xi)&=\sqrt{\frac{1}{\varphi(\xi)}}m_{\ell}(\xi),\   1 \le \ell \le qN-1.
\end{align*}

Then, we have
\begin{align}
\widehat{\widetilde{\psi_0}}(\mathfrak{p}^{-1}N\xi)&=\sqrt{\varphi(\mathfrak{p}^{-1}N\xi)}\hat{\psi_0}(\mathfrak{p}^{-1}N\xi)\notag \\\
&=\sqrt{\varphi(\mathfrak{p}^{-1}N\xi)}\{m_0(\xi)\hat{\psi_0}(\xi)\}\notag\\\
&=\sqrt{\varphi(\mathfrak{p}^{-1}N\xi)}\left\{m_0(\xi)\frac{\widehat{\widetilde{\psi_0}}(\xi)}{\sqrt{\varphi(\xi)}}\right\}\notag\\\
&=\sqrt{\frac{\varphi(\mathfrak{p}^{-1}N\xi)}{\varphi(\xi)}}m_0(\xi)\widehat{\widetilde{\psi_0}}(\xi)\notag\\\
&=\widetilde{m_0}(\xi)\widehat{\widetilde{\psi_0}}(\xi) \tag{3.19},
\end{align}
and

\begin{align}
\lim\limits_{\xi \to 0^{+}} \widehat{\widetilde{\psi_0}}(\xi)=\lim\limits_{\xi \to 0^{+}}\sqrt{\varphi(\xi)}\hat{\psi_0}(\xi)= 1.\tag{3.20}
\end{align}
Since $\{\psi_{\ell}, m_{\ell}\}_{\ell=0}^{qN-1}$ is a general setup, by (3.18), we have
\begin{align}
\text{Supp}\   \widehat{\widetilde{\psi_0}}(\xi) \subseteq q^2N\mathfrak{D}, \tag{3.21}
\end{align}
and

\begin{align}
\sum\limits_{\ell=0}^{qN-1}|\widetilde{m_{\ell}}(\xi)|^2&=|\widetilde{m_0}(\xi)|^2+\sum\limits_{\ell=1}^{qN-1}|\widetilde{m_{\ell}}(\xi)|^2\notag
\\
&=\frac{\varphi(\mathfrak{p}^{-1}N\xi)}{\varphi(\xi)}|m_0(\xi)|^2+\sum\limits_{\ell=1}^{qN-1}\frac{|m_{\ell}(\xi)|^2}{\varphi(\xi)}\notag\\
&=\frac{1} {\varphi(\xi)}\varphi(\xi)\notag\\
&=1.\tag{3.22}
\end{align}
Thus

\begin{align}
\widetilde{m_{\ell}}(\xi)\in L^{\infty}(\mathbb{K}) \ \text{for} \   1 \le \ell \le qN-1. \tag{3.23}
\end{align}
Let  $\widetilde{\psi}_1, \widetilde{\psi}_2,\dots, \widetilde{\psi}_{qN-1} \in L^2(\mathbb{K})$ be such that

\begin{align}
\widehat{\widetilde{\psi}}_{\ell}(\mathfrak{p}^{-1}N\xi)= \widetilde{m}_{\ell} (\xi)\widehat{\widetilde{\psi_0}}(\xi),  1 \le \ell \le qN-1.\tag{3.24}
\end{align}

Define
\begin{align*}
\widetilde{\cal M}(\xi)=\begin{bmatrix}
\widetilde{{m_0}}(\xi)\\
\widetilde{{m_1}}(\xi)\\
\vdots\\
\widetilde{m}_{qN-1}(\xi)
\end{bmatrix}_{qN \times 1}.
\end{align*}
Then,  by (3.19), (3.20), (3.21) and (3.23),  the collection
 $\{\widetilde{\psi_{\ell}}, \widetilde{m_{\ell}}\}_{\ell=0}^{qN-1}$ is  a general setup.

Using (3.22), we have

\begin{align*}
\widetilde{\cal M}(\xi)^\ast \widetilde{\cal M}(\xi)
&=\sum\limits_{\ell=0}^{qN-1}|\widetilde{m_{\ell}}(\xi)|^2=1.
\end{align*}
Hence, by Theorem 3.1, $\{D_{\mathfrak p^{j}}T_{\lambda}\widetilde{\psi}_{\ell}\}_{j\in \mathbb{Z},\lambda \in \Lambda \atop 1\le \ell \le qN-1}$ is a Parseval nonuniform wavelet frames for $L^2(\mathbb{K})$.

Next, we compute
\begin{align*}
\hat{\psi_{\ell}}(\mathfrak{p}^{-1}N\xi)&=m_{\ell}(\xi)\hat{\psi_0}(\xi)\\
&=\left\{\widetilde{m_{\ell}}(\xi)\sqrt{\varphi(\xi)}\right\}\left\{\frac{\widehat{\widetilde{\psi_0}}(\xi)}{\sqrt{\varphi(\xi)}}\right\}\\
&=\widetilde{m_{\ell}}(\xi)\widehat{\widetilde{\psi_0}}(\xi)\\
&=\widehat{\widetilde{\psi_{\ell}}}(\mathfrak{p}^{-1}N\xi).
\end{align*}
This gives,  $\psi_{\ell}=\widetilde{\psi_{\ell}}$. Hence, the system  $\{D_{\mathfrak p^{j}}T_{\lambda}\psi_{\ell}\}_{j\in \mathbb{Z},\lambda \in \Lambda \atop  1 \le \ell \le qN-1}$ is a Parseval nonuniform wavelet frames for $L^2(\mathbb{K})$.
\end{proof}

{\bf {Example 3.1}} Let $\psi(x)\in L^2(\mathbb{K})$ such that $\hat{\psi_0}(\xi)={\bf 1}_{q^2N\mathfrak{D}}(\xi) $,then

$$\text{Supp} \ \hat{\psi_0}(\xi) \subseteq q^2N\mathfrak{D},~\lim \limits_{\xi \to 0^+} \hat{\psi_0}(\xi)=1,$$
$$\hat{\psi_0}(\mathfrak{p}^{-1}N\xi) ={\bf 1}_{q^2N\mathfrak{D}}(\mathfrak{p}^{-1}N\xi)={\bf 1}_{q^3N^2\mathfrak{D}} (\xi)={\bf 1}_{q^3N^2\mathfrak{D}} (\xi) (\xi){\bf 1}_{q^2 N\mathfrak{D}} (\xi) = m_0(\xi)\hat{\psi_0}(\xi),$$
where $m_0(\xi)={\bf 1}_{q^3N^2\mathfrak{D}} (\xi) \in L^{\infty}(\mathbb{K})$.

Define
\begin{align*}
m_1(\xi)={\bf 1}_{\mathbb{K}-q^3 N^2\mathfrak{D}}(\xi), \quad m_2(\xi)=1 \  \text{and} \  \quad m_3(\xi)=i.
\end{align*}
We choose $\psi_1,\psi_2,\psi_3 \in L^2({\mathbb{K}})$ and
\begin{align*} {\cal M}(\xi)=\begin{bmatrix}
m_0(\xi)\\
m_1(\xi)\\
m_2(\xi)\\
m_3(\xi)
\end{bmatrix} 
\end{align*}
such that
$$  \hat{\psi_{\ell}}(\mathfrak{p}^{-1}N\xi)=m_{\ell}(\xi)\hat{\psi_0}(\xi), \ \ell=1, 2, 3.$$
Then, the collection  $\{\psi_\ell,m_\ell\}_{\ell=0}^3$ is a general setup such that
\begin{align*}
{{\mathcal M}(\xi)}^\ast {\cal M}(\xi)=\displaystyle\sum_{\ell=0}^{3}|m_\ell(\xi)|^2=1.
\end{align*}
 Hence, by Theorem 3.1, the family $\{{D_{\mathfrak{p}^j}}T_{\lambda}\psi_{\ell}\}_{j\in \mathbb{Z},\lambda \in \Lambda \atop \ell=1,2,3}$ is a nonuniform Parseval  wavelet frame for $L^{2}(\mathbb{K})$.\\
 
 \parindent=0mm \vspace{.1in}
 {\bf{4. Potential Applications}}

\parindent=0mm \vspace{.1in}
The results specified in this article are theoretical in nature and will definitely promote new directions to the modern theory of Wavelet analysis and broaden its field of applications. The presented work will be of substantial importance to that part of physical and mathematical community working in the theory of harmonic analysis, chaotic systems, quasi-crystal theory, signal and image processing, data transmission with erasures, quantum computing, medicine, representation theory and algebraic geometry.

\parindent=0mm \vspace{.1in}
{{(a) \it Chaotic Systems}}: We deal with non-Archimedian  fields, that is, the norm satisfies a stronger inequality than the triangle inequality namely \textit{ultrametric inequality}. The $p$-adic distance leads to interesting deviations from the classical real analysis, the geometry of these spaces is unlike the euclidian geometry based on real space $\mathbb R$. In non-Archimedian geometry two different balls are either disjoint or the one is contained in the other one (splitting property). Furthermore the field of 2-adic and 2-series numbers have a hierarchical structure: every disc consists of two disjoint discs of smaller radius (tree property). Thus these fields are homeomorphic to a Cantor set on $\mathbb R$. The fractal-like structure of these fields enable their application not only for the description of geometry at small distances, but also for describing chaotic behavior of chaotic
systems.

\parindent=0mm \vspace{.1in}
{(b) \it Quasi-crystal Theory}: The description of the elements of a vector space based on the use of an overcomplete
system is a general method rediscovered several times in different areas of mathematics, science
and engineering. For example, in crystallography there exists an alternative description for
the hexagonal crystals based on the use of an additional axis. The use of a frame
leads to a simpler description of atomic positions in a diamond type crystal. This leads to
a simpler description of the symmetry transformations and of the mathematical objects with
physical meaning. Some of the most important models used in quasicrystal physics can be
generated in a unitary way by using the imbedding into a superpace defined by certain frames.
These observations allow a fructuous interchange of ideas and methods between frame theory
and quasicrystal physics.

\parindent=0mm \vspace{.1in}
{(c) \it Application to Geophysics}:   $p$-adics can be applied  to geophysics using a p-adic diffusion representation of the master equations for the dynamics of a fluid in capillaries in porous media and formulate several mathematical problems motivated by such applications.  $p$-adic wavelets are a powerful tool for obtaining analytic solutions of diffusion equations. Because $p$-adic diffusion is a special case of fractional diffusion, which is closely related to the fractal structure of the configuration space, $p$-adic geophysics can be regarded as a new approach to fractal modeling of geophysical processes.

\parindent=0mm \vspace{.2in}

{\bf{References}}

\end{document}